\documentstyle[12pt]{amsart}
\newtheorem{theorem}{Theorem}[section]
\newtheorem{corollary}[theorem]{Corollary}
\newtheorem{remark}[theorem]{Remark}

\newtheorem{proposition}[theorem]{Proposition}
\newtheorem{definition}{Definition}[section]

\numberwithin{equation}{section}

\begin{document}

\title [Symplectic and Lagrangian mean curvature flows]
 {Singularities of symplectic and Lagrangian mean curvature flows}
\author{Xiaoli Han, Jiayu Li}

\address{Math. Group, The abdus salam ICTP\\ Trieste 34100,
   Italy}
\email{xhan@@ictp.it}

\address{Math. Group, The abdus salam ICTP\\ Trieste 34100,
   Italy\\
   and Academy of Mathematics and Systems Sciences\\ Chinese Academy of
Sciences\\ Beijing 100080, P. R. of China. } \email{jyli@@ictp.it}


\date{}

\begin{abstract}
In this paper we study the singularities of the mean curvature flow
from a symplectic surface or from a Lagrangian surface in a
K\"ahler-Einstein surface. We prove that the blow-up flow
$\Sigma_s^\infty$ at a singular point $(X_0, T_0)$ of a symplectic
mean curvature flow $\Sigma_t$ or of a Lagrangian mean curvature
flow $\Sigma_t$ is a non trivial minimal surface in ${\bf R}^4$, if
$\Sigma_{-\infty}^\infty$ is connected.
\end{abstract}

\maketitle

\section{Introduction}
Suppose that $M$ is a compact K\"ahler-Einstein surface. Let
$\omega$ be the K\"ahler form on $M$ and $\langle\cdot,
\cdot\rangle$ be the K\"ahler metric, the K\"ahler angle $\alpha$
of $\Sigma$ in $M$ is defined by
\begin{equation}\label{0.1}\omega|_\Sigma=\cos\alpha
d\mu_\Sigma\end{equation} where $d\mu_\Sigma$ is the area element
of $\Sigma$ of the induced metric from $\langle,\rangle$. We call
$\Sigma$ a symplectic surface if $\cos\alpha>0$, $\Sigma$ a
Lagrangian surface if $\cos\alpha =0$, and call $\Sigma$ a
holomorphic cure if $\cos\alpha\equiv 1$.

It is proved in \cite{CT}, \cite{CL1} and \cite{Wa} that, if the
initial surface is symplectic, then along the mean curvature flow,
at each time $t$ the surface is still symplectic, which we call a
symplectic mean curvature flow. It is proved in \cite{Sm1} that,
if the initial surface is Lagrangian, then along the mean
curvature flow, at each time $t$ the surface is still Lagrangian,
which we call a Lagrangian mean curvature flow.

We \cite{HL} showed that, if the scalar curvature of the
K\"ahler-Einstein surface is positive and the initial surface is
sufficiently close to a holomorphic curve, then the mean curvature
flow has a global solution and it converges to a holomorphic
curve.

In general, the mean curvature flows may produce singularities.
The beautiful results on the nature of singularities of the mean
curvature flows of convex hypersurfaces have been obtained by
Huisken \cite{H2}, Huisken-Sinestrari \cite{HS1},\cite{HS2} and
White \cite{W}. In \cite{HL1} we obtain the relation between the
maximum of the K\"ahler anlge and the maximum of $|H|^2$ on the
blow-up flow of the symplectic mean curvature flow or the
calibrated Lagrangian mean curvature flow.

It is well-known (see \cite{HL1}) that, a sequence of rescaled
surfaces at a singular point $(X_0, T_0)$ converges strongly to a
blow-up flow $\Sigma^\infty_s$, for $s\in (-\infty, 0]$ if the
singular point is of type I, for $s\in (-\infty, +\infty)$ if the
singular point is of type II. In this paper, we show that
$\Sigma_s^\infty$ is a non trivial holomorphic curve with finite
total curvature and bounded Gauss curvature in ${\bf C}^2$, if
$\Sigma_{-\infty}^\infty$ is connected. It is well known that the
total curvature of a non-flat minimal surface with finite total
curvature in ${\bf C}^2$ achieves only discrete values $-2\pi N$
where $N$ is a nature number. We therefore believe that the size
of the singular set can be controlled.

More precisely, let $T$ be a discrete singular time and $(X_0, T)$
be a singular point in $M$, one shows (see \cite{HL1}) that there
are sequences $r_k\to 0$, $0<\sigma_k\leq r_k/2$, $t_k\in
[T-(r_k-\sigma_k)^2, T- r_k^2/4]$,
$F(x_k,t_k)=X_k\in\bar{B}_{r_k-\sigma_k}(X_0)$, such that
$$\lambda_k^2=|A|^2(X_k)=|A|^2(x_k, t_k)=\sup_{[T-(r_k-\sigma_k)^2, T- r_k^2/4]}
\sup_{\Sigma_t\cap B_{r_k-\sigma_k}(X_0)}|A|^2,$$ and
$$\sup_{[t_k-(\sigma_k/2)^2, t_k]}
\sup_{\Sigma_t\cap B_{r_k-\sigma_k/2}(X_0)}|A|^2\leq 4\lambda_k^2.
$$
Choose a normal coordinates in a neighborhood of $X_0$, express
$F$ in this coordinates, and consider the following sequence of
rescaled surfaces
$$F_k(x, s)=\lambda_k(F(x_k+x, t_k+\lambda_k^{-2}s)-F(x_k, t_k)),~~~~~~
s\in [-\lambda_k^2\sigma_k^2/4, 0].$$ It is clear that
$|A_k|(0,0)=1$ and $|A_k|^2\leq 4$. Denote the rescaled surfaces
by $\Sigma^k_s$, then $\Sigma^k_s\to\Sigma^\infty_s$ in
$C^2(B_R(X_0)\times [-R, R])$ for any $R>0$ and any
$B_R(X_0)\subset {\bf R}^4$. We call $\Sigma^\infty_s$ a blow-up
flow.

If there exists $R_0>0$ such that, for all $R>R_0$,
$\Sigma_s^{\infty}\cap B_R(0)\not =\emptyset$ for sufficiently
large $-s$, then $\Sigma_s^{\infty}$ converges in $C^2$ to
$\Sigma_{-\infty}^{\infty}$ in $B_R(0)$ as $s\to\infty$ for all
$R>R_0$.

We prove the following main theorem. In fact, we prove a more
general result (see Theorem \ref{maintheorem}) which implies the
following one.

\vspace{.2in}

\noindent {\bf Main Theorem} {\it Let $M$ be a K\"ahler-Einstein
surface and $\Sigma_0$ be a symplectic surface in $M$. If the limit
$\Sigma_{-\infty}^{\infty}$ of the blow-up flow $\Sigma^\infty_s$ at
infinity is connected, then
 $\Sigma^\infty_s$ is independent of
$s$ denoted by $\Sigma^\infty$, and $\Sigma^\infty$ is a non
trivial holomorphic curve in ${\bf C}^2$ with Gauss curvature
$-2\leq K\leq 0$, and finite total curvature
$$
-\int_{\Sigma^\infty}KdV=2\pi N,
$$
where $N$ is a nature number.}

\vspace{.2in}

\noindent {\bf Remark} {\it Let $A_\infty$ be the second
fundamental form of $\Sigma^\infty$ in ${\bf C}^2$, then
$|A_\infty|^2=-2K$, so we can also write the last identity as
$$
\int_{\Sigma^\infty}|A_\infty|^2dV=4\pi N.
$$
}

\vspace{.2in}

Due to the theorem, we believe that, once we control
$\int_{\Sigma_t}|A|^2d\mu_t$, we can control the singular size of
the mean curvature flow. Based on it, we propose a conjecture at
end.

\vspace{.2in}

\noindent {\bf Conjecture} {\it A symplectic mean curvature flow in a
K\"ahler-Einstein surface blows up at most countable discrete times and
at each blow-up time, the blow-up set consists of at most finitely
many points. }

\vspace{.2in}

We prove some similar results in the case of the Lagrangian mean
curvature flows.

\section{Properties of blow-up flows of symplectic mean curvature flows}

In this section, we prove our main theorem.

Let $T$ be an isolated singular time, that is, the mean curvature
flow exists in $t\in [T-\epsilon, T)$, and $(X_0, T)$ be a blow-up
point. From the main theorem in \cite{CL1} and \cite{Wa}, we know
that this is a type II singularity. Recall that \cite{HL1}, we can
define a sequence of rescaled surfaces $\Sigma^k_s$ around
$(X_0,T)$. For each fixed $R>0$, by parabolic estimates, we have
that $\Sigma^k_s\to\Sigma^\infty_s$ in $C^2(B_R(0)\times [-R,R])$
for any $B_R(0)\subset {\bf C}^2$, and $\Sigma_s^\infty$ also
evolves along the mean curvature flow with the property that
\begin{equation}\label{abound}
|A|^2(0,0)=1,~~{\rm and}~~|A|^2\leq 4. \end{equation} By the
evolution equation derived in \cite{CL1}, we see that, along the
mean curvature flow $\Sigma_s^\infty$, $\cos\alpha$ satisfies
$$(\frac{\partial}{\partial t}-\Delta)\cos\alpha=
|\overline{\nabla}_0J_{\Sigma^\infty_s}|^2\cos\alpha ,$$ where
$\overline{\nabla}_0$ is the classical derivative in ${\bf R}^4$.

By the monotonicity inequality (Proposition 3.2 in \cite{HL1}), we
have, for any $R>0$, $s<0$,
\begin{equation}\label{e0.0} \mu_s^\infty(\Sigma_s^\infty\cap
B_R(0))\leq CR^2,
\end{equation}
where $C>0$ does not depend on $s$ or $R$. Since $M$ is compact,
then there exists a constant $\delta$ such that $\cos\alpha\geq
\varepsilon_0$ on $M$. It is easy to see that $\cos\alpha$ is
scaling invariant, thus $\cos\alpha\geq \varepsilon_0$ on
$\Sigma^\infty_s$, for all $s\in (-\infty, +\infty)$. Therefore,
on $\Sigma^\infty_s$, the Isoperimetric inequality holds. That is,
\begin{proposition}
There is a positive constant $C(\epsilon_0)$ which depends only on
$\epsilon_0$ such that, for any open smooth (connected) domain
$A\subset \Sigma^\infty_s$,
\begin{equation}\label{2.8}
{\rm Area }(A)\leq C(\epsilon_0) ({\rm Length}( \partial A))^2.
\end{equation}
\end{proposition}

{\it Proof.} By Theorem 30.1 in [Si1], there is an integral
current $B$ with compact support such that $\partial B=\partial A$
and
$$
{\rm Area }(B)\leq C ({\rm Length}( \partial A))^2,
$$
where $C$ is an absolute positive constant.

Let $T$ be the cone over $A-B$ with $\partial T=A-B$. Since
$d\omega =0$ and $\omega |_{\Sigma^\infty_s}=\cos\alpha
d\mu_{s}^\infty$, we have
\begin{eqnarray*}
{\rm Area }(A)&\leq &\frac{1}{\epsilon_0}\int_A\omega\\
&=&\frac{1}{\epsilon_0}\int_B\omega +\partial T(\omega)
=\frac{1}{\epsilon_0}\int_B\omega \\
&\leq & \frac{1}{\epsilon_0}{\rm Area }(B)\\
&\leq & C(\epsilon_0) ({\rm Length}( \partial A))^2.
\end{eqnarray*}
\hfill Q. E. D.

Fix $R>0$. For any point $x$ in the connected components of
$\Sigma^\infty_s\cap B_R(0)$ that intersect with $B_{R/2}(0)$,
denote the intrinsic ball of radius $r$ around $x$ by $\hat
B_r(x)$. The isoperimetric inequality implies that
$${\rm Vol} (\hat B_r(x))\geq C r^2,$$ where $C$ is a constant which
depends only on $\varepsilon_0$. By (\ref{e0.0}), we see that
$\Sigma_s^\infty\cap B_R(0)$ contains at most finite many
connected components which intersect with $B_{R/2}(0)$. We denote
it by $\Sigma_s^{(\infty,l)}, l=1,\cdots, L$.

\begin{proposition}
If the blow-up flow of a symplectic mean curvature flow is minimal,
that is $H\equiv 0$, it must be holomorphic.
\end{proposition}
{\it Proof:} From \cite{CW} we know that on the minimal surface
$\cos\alpha$ satisfies that,
\begin{equation}\label{e1}
-\Delta\cos\alpha=2|\nabla\alpha|^2\cos\alpha=2\frac{|\nabla\cos\alpha|^2}{1-\cos^2\alpha}
\cos\alpha, \end{equation} at the points which are not
holomorphic. Fix $R>0$, let $\phi(x)\in C^\infty(B_{2R}(0))$ be a
cut-off function such that $\phi\equiv 1$ in $B_R(0)$ and
$\phi\equiv 0$ outside of $B_{2R}(0)$. Multiplying the equation
(\ref{e1}) by $(1-\cos^2\alpha)\phi^2$ and integrating by parts,
we get that,
\begin{eqnarray}\label{e2}
2\int
|\nabla\cos\alpha|^2\cos\alpha\phi^2=\int\sin^2\alpha\phi\nabla\phi\cdot\nabla\cos\alpha.
\end{eqnarray}
By Schwartz inequality,
\begin{eqnarray}\label{e3}
\int\sin^2\alpha\phi\nabla\phi\cdot\nabla\cos\alpha\leq
\frac{1}{2}\int\frac{\sin^4\alpha}{\cos\alpha}|\nabla\phi|^2+\frac{1}{2}\int|\nabla\cos\alpha|^2
\cos\alpha\phi^2.
\end{eqnarray}
Plugging inequality (\ref{e3}) into (\ref{e2}), we get,
\begin{eqnarray*}
\int|\nabla\cos\alpha|^2\cos\alpha\phi^2\leq\int\frac{\sin^4\alpha}{\cos\alpha}|\nabla\phi|^2.
\end{eqnarray*}
(\ref{e0.0}) yields,
\begin{eqnarray*}
\int_{\Sigma^\infty\cap B_R(0)}|\nabla\cos\alpha|^2\leq
C(\varepsilon_0)\frac{vol(B_{2R})}{R^2}\leq C.
\end{eqnarray*}
Let $R\to\infty$, we get that,
\begin{eqnarray}\label{e4}
\int_{\Sigma^\infty}|\nabla\cos\alpha|^2\leq C.
\end{eqnarray}
Multiplying equation (\ref{e1}) by $\cos^p\alpha\phi^2$, where
$p>0$ will be determined later.
\begin{eqnarray*}
2\int\cos^{p+1}\alpha|\nabla\alpha|^2\phi^2-p\int\cos^{p-1}\alpha|\nabla\cos\alpha|^2\phi^2
=2\int\cos^p\alpha\phi\nabla\phi\cdot\nabla\cos\alpha.
\end{eqnarray*}
Using Holder inequality , (\ref{e4}), and (\ref{e0.0}), we have
\begin{eqnarray*}
2\int\cos^p\alpha\phi\nabla\phi\nabla\cos\alpha&\leq& 2(\int
\cos^{2p}\alpha|\nabla\phi|^2\phi^2)^{1/2}(\int_{B_{2R}/B_R}|\nabla\cos\alpha|^2)^{1/2}\\
&\leq& C(\int_{B_{2R}/B_R}|\nabla\cos\alpha|^2)^{1/2}\to
0,~~~~~~~{\rm as}~~~~~~~~R\to\infty.
\end{eqnarray*} Thus we have,
\begin{eqnarray*}
\int(2\cos^2\alpha-p\sin^2\alpha)|\nabla\alpha|^2\cos^{p-1}\alpha=0.
\end{eqnarray*}
Choosing $p$ such that $\frac{p}{p+2}<\varepsilon_0^2$, then
$2\cos^2\alpha-p\sin^2\alpha>c(\varepsilon_0)>0$, thus
\begin{eqnarray*}
\int |\nabla\alpha|^2\cos^{p-1}\alpha=0,
\end{eqnarray*} which implies that $|\nabla\alpha|^2=0$.
Therefore, $\Sigma$ is holomorphic with some complex structure in
${\bf C}^2$.

\hfill  Q. E. D.

\begin{remark} This result can also be deduced from the theorem in
\cite{CY}.
\end{remark}

\vspace{.2in}

 It is clear that, if
$\Sigma_s^{(\infty,l)}\cap B_R(0)\not =\emptyset$ for all $-s$
sufficiently large, then $\Sigma_s^{(\infty,l)}$ converges in
$C^2(B_R(0))$ to $\Sigma_{-\infty}^{(\infty,l)}$ as $s\to-\infty$,
for any $R>0$.

\begin{definition}
The component $\Sigma_s^{(\infty,l)}$ is called {\it simple}, if
$\Sigma_{-\infty}^{(\infty,l)}$ is connected.
\end{definition}

\begin{theorem}\label{maintheorem} Let $M$ be a K\"ahler-Einstein
surface and $\Sigma_0$ be a symplectic surface in $M$. Each simple
connected component $\Sigma_{s}^{(\infty,l)}$ of the blow-up flow
$\Sigma^\infty_s$ is independent of $s$ denoted by
$\Sigma^{(\infty,l)}$, and $\Sigma^{(\infty,l)}$ is a holomorphic
curve in ${\bf C}^2$ with Gauss curvature $-2\leq K\leq 0$, and
finite total curvature
$$
-\int_{\Sigma^\infty}KdV=2\pi N,
$$
where $N$ is a nature number.
\end{theorem}

{\it Proof:} For simplicity, we denote the simple connected
component by $\Sigma^\infty_s$. The following monotonicity formula
for $\Sigma^\infty_s$ is derived in \cite{HL1}.

Let $H (X, X_0, t, t_0)$ be the backward heat kernel on ${\bf
R}^4$. Define
$$
\rho (X, X_0, t, t_0)=4\pi (t_0-t)H (X, X_0, t, t_0)=\frac{1}{4\pi
(t_0-t)}\exp \left(-\frac{|X-X_0|^2} {4(t_0-t)}\right)
$$
for $t<t_0$. We have, for $-\infty<s<s_0$,
\begin{eqnarray}\label{mon4}
&&\frac{\partial}{\partial
s}\left(\int_{\Sigma^\infty_s}\frac{1}{\cos\alpha} \rho(F_\infty,
0, s, s_0)d\mu^\infty_s\right)\nonumber\\&&= -\left(
\int_{\Sigma^\infty_s}\frac{1}{\cos\alpha}\phi\rho (F_\infty, 0,
s, s_0)
\left|H_\infty+\frac{(F_\infty)^{\perp}}{2(s_0-s)}\right|^2d\mu_s^\infty \right.\nonumber \\
&&\left.+ \int_{\Sigma^\infty_s}\frac{1}{2\cos\alpha}\phi\rho
(F_\infty,0,s,s_0)\left|\overline{\nabla}J_{\Sigma^\infty_s}\right|^2d\mu_s^\infty
\right.\nonumber\\
&&\left.+\int_{\Sigma^\infty_s}\frac{2}{\cos^3\alpha}\left|\nabla
\cos\alpha\right|^2\phi\rho
(F_\infty,0,s,s_0)d\mu_s^\infty\right).
\end{eqnarray}

Choosing $s_0=0$, $s_1=4t$, $s_2=2t$, for $t<0$ we have,
\begin{eqnarray}\label{bb}
&&\int_{\Sigma^\infty_{4t}}\frac{1}{\cos(x, 4t)}
\frac{1}{-4t}e^{-\frac{|F_\infty|^2}{-4t}}d\mu^\infty_{4t}-\int_{\Sigma^\infty_{2t}}
\frac{1}{\cos(x, 2t)}
\frac{1}{-2t}e^{-\frac{|F_\infty|^2}{-2t}}d\mu^\infty_{2t}\nonumber\\&&\geq
\int_{4t}^{2t}\int_{\Sigma^\infty_s}\frac{1}{\cos\alpha}\rho(F_\infty,0,s,0)
\left|\overline{\nabla}J_{\Sigma^\infty_s}\right|^2d\mu^\infty_s
ds.\nonumber
\end{eqnarray}
Since $\int_{\Sigma^\infty_s}\frac{1}{\cos\alpha} \rho(F_\infty,
0, s, s_0)d\mu^\infty_s$ is uniformly bounded above (See
\cite{HL1}), so the left side of the above inequality tends to
zero as $t\to -\infty$. Moreover,
\begin{eqnarray*}
&&\int^{2t}_{4t}\int_{\Sigma_s^\infty}\frac{1}{\cos\alpha}\frac{1}{-s}
e^{-\frac{|F_\infty|^2}{-s}}|\overline\nabla J_{\Sigma^\infty_s}|^2 d\mu_s^\infty\\
&&=-2t\int_{\Sigma_{s'}^\infty}\frac{1}{\cos\alpha}\frac{1}{-s'}
e^{-\frac{|F_\infty|^2}{-s'}}|\overline\nabla J_{\Sigma^\infty_{s'}}|^2 d\mu_{s'}^\infty\\
&&\geq C\int_{\Sigma_{s'}^\infty} |\overline\nabla
J_{\Sigma^\infty_{s'}}|^2
e^{-\frac{|F_\infty|^2}{-s'}}d\mu_{s'}^\infty
\end{eqnarray*}
where $s'\in [4t, 2t]$, $C$ is independent of $t$.

It is clear that, for any $R>0$, if $\Sigma^\infty_s\cap
B_R(0)\neq\emptyset$ for sufficiently large $-s$, then
$\Sigma_s^\infty\cap B_R(0)$ converges strongly to
$\Sigma_{-\infty}^\infty\cap B_R(0)$ as $s\to -\infty$. By the
assumption, we know that $\Sigma_{-\infty}^\infty\cap B_R(0)$ is
connected for $R$ large enough, so there is $S_0$ such that
$s<S_0$, $\Sigma_{s}^\infty\cap B_R(0)$ is connected.

Letting $t\to -\infty$, we get that on $\Sigma^\infty_{-\infty}$,
$$\left|\overline{\nabla}J_{\Sigma_{-\infty}^\infty}\right|^2=0,\,\,\,\,\,
{\rm that\,is,}\,\,\, |\nabla \cos\alpha|^2=0.$$ It follows that
$\cos\alpha\equiv\theta_0$, $\theta_0$ is constant on
$\Sigma^\infty_{-\infty}$.

We can choose the suitable complex structure on ${\bf R}^4$ such
that with respect to this complex structure we have
$\cos\alpha\equiv 1$ on $\Sigma^\infty_{-\infty}$. In fact, assume
that $\omega$ is written as $\omega=dz_1\wedge d\bar
z_1+dz_2\wedge d\bar z_2$ under the standard complex structure of
${\bf R}^4$. We define a new complex structure $J^\ast$ of ${\bf
R}^4$ as follows:
$$J^\ast(\partial/\partial x_1)=\theta_0(\partial/\partial y_1),
~~~~~~J^\ast(\partial/\partial y_1)=-1/\theta_0(\partial/\partial
x_1),$$
$$J^\ast(\partial/\partial x_2)=1/\theta_0(\partial/\partial y_2),
~~~~~~J^\ast(\partial/\partial y_2)=-\theta_0(\partial/\partial
x_2).$$ Under the complex structure $J^\ast$, the complex
coordinates are $z_1^\ast=x_1+\sqrt{-1}\theta_0^{-1}y_1$,
$z_2^\ast=\theta_0^{-1}x_2+\sqrt{-1}y_2$. Thus $\omega^\ast=
dz_1^\ast\wedge d\bar z^\ast_1+dz_2^\ast\wedge d\bar z_2^\ast$
satisfies
$\omega^\ast|_{\Sigma^\infty_{-\infty,l}}=d\mu_{-\infty,l}^\infty$,
In other words, $\theta_0\equiv 1$.

Recall that, on $\Sigma^\infty_s$, $\cos\alpha$ satisfies
$$(\frac{\partial}{\partial s}-\Delta)\cos\alpha=
|\overline{\nabla}_0J_{\Sigma^\infty_s}|^2\cos\alpha.$$ Using the
maximum principle on $\Sigma^\infty_s$ (See \cite{EH}), for all
$s\leq S_0$, we can see that $\min_{\Sigma^\infty_{s}}\cos\alpha$
is a nondecreasing function of $s$, so $\cos\alpha\equiv 1$ on
$\Sigma_{s}^\infty$, for all $s\leq S_0$. We therefore have that
$|\overline{\nabla}_0J_{\Sigma^\infty_s}|\equiv 0$, hence
$H_s^\infty\equiv 0$, and consequently, $\Sigma^\infty_s$, $s\leq
S_0$, does not depend on $s$ and is a holomorphic curve. We
continue this process and claim that $\Sigma^\infty_s$ is a
holomorphic curve for all $s<0$.

By (\ref{abound}), we see that the second fundamental form
$A_\infty$ of $\Sigma^\infty$ in ${\bf R}^4$ satisfies
$$
|A_\infty(0)|=1~~{\rm and}~~|A_\infty|^2\leq 4.
$$

Let $K$ be the Gauss curvature of $\Sigma^\infty$, $R$ the
curvature operator of $M$. By Gauss equation,
$$
K_{1212}=R_{1212}+(h_{11}^\alpha h_{22}^\alpha - h_{12}^\alpha
h_{12}^\alpha),
$$
we get
$$
|A|^2=|H|^2-2K_{1212}+2R_{1212}.
$$

Thus we have
$$
K=-\frac{1}{2}|A_\infty|^2.
$$
We will show in the following proposition that, the total
curvature of $\Sigma^\infty$ is finite.

\begin{proposition}\label{pro1}
The minimal surface $\Sigma^\infty$ in the previous theorem is of
finite type, that is, its total curvature is finite.
\end{proposition}
{\it Proof.} Since
$$-\int_{\Sigma^\infty}K d\mu^\infty=\frac{1}{2}\int_{\Sigma^\infty}|A_\infty|^2
d\mu^\infty,
$$
so it suffices to prove that
$$\int_{\Sigma^\infty}|A_\infty|^2
d\mu^\infty<\infty .
$$

It follows from the integral curvature estimate (Theorem $4$, in
[Il]) that for any $r>0$,
$$r^{-2}\int_{t-r^2}^t\int_{\Sigma_t\cap B_r(x)} |A|^2 d\mu_t dt\leq
C.
$$
Here $C$ does not depend on $r$. Thus as $k$ sufficiently large we
have
$$\lambda_k^{2}R^{-2}\int_{t_k-\lambda_k^{-2}R^2}^{t_k}\int_{\Sigma_t\cap B_{\lambda_k^{-1}R}(X_k)}
|A|^2 d\mu_t dt\leq C.
$$
This implies that
$$R^{-2}\int_{-R^2}^{0}\int_{\Sigma^k_s\cap B_{R}(0)}
|A_k|^2 d\mu^k_s ds\leq C.
$$
Letting $k\to\infty$, we get
$$R^{-2}\int_{-R^2}^{0}\int_{\Sigma^\infty_s\cap B_{R}(0)}
|A_\infty|^2 d\mu^\infty_s ds\leq C.
$$
Note that $\Sigma^\infty_s$ does not depend on $s$, we have
$$\int_{\Sigma^\infty\cap B_{R}(0)}
|A_\infty|^2 d\mu^\infty \leq C.
$$
Letting $R\to\infty$, we get that
$$\int_{\Sigma^\infty}
|A_\infty|^2 d\mu^\infty \leq C.
$$
This proves the proposition. \hfill Q. E. D.

By Proposition 6.1 in \cite{HO}, we know that
$$
-\int_{\Sigma^\infty}KdV=2\pi N.
$$
where $N$ is a natural number.

This proves the theorem. \hfill Q. E. D.

\begin{corollary} Let $M$ be a K\"ahler-Einstein
surface and $\Sigma_0$ be a symplectic surface in $M$. If the
blow-up flow $\Sigma^\infty_s$ contains only one connected
component, and its limit at $-\infty$, $\Sigma_{-\infty}^\infty$ is
connected, then
 $\Sigma^\infty_s$ is independent of
$s$ denoted by $\Sigma^\infty$, and $\Sigma^\infty$ is a non trivial
holomorphic curve in ${\bf C}^2$ with Gauss curvature $-2\leq K\leq
0$, and finite total curvature
$$
-\int_{\Sigma^\infty}KdV=2\pi N,
$$
where $N$ is a nature number.
\end{corollary}

There are many works on minimal surfaces with finite total
curvature, we state one of them (c.f. Theorem 6.1 in \cite{O}), so
that the readers are aware of the properties of $\Sigma^\infty$.
\begin{theorem}
Let $x: \Sigma\to {\mathbb R}^m$ be a non-flat complete minimal
surface in ${\mathbb R}^m$. Then, the following conditions are
equivalent.

(i) $M$ has finite total curvature.

(ii) $M$ is biholomorphic with a compact Riemann surface $\bar{M}$
with finitely many points removed and each $\omega_i=\partial x_i$
extends to $\bar{M}$ as a meromorphic form.

(iii) $M$ is biholomorphic with an open subset of a compact
Riemann surface $\bar{M}$ and the Gauss map $G: M\to P^{m-1}(C)$
extends to a holomorphic map of $\bar{M}$ into $P^{m-1}(C)$.
\end{theorem}

It is certainly important to know when a connected component
$\Sigma_s^{(\infty,l)}$ of the blow-up flow $\Sigma_s^\infty$ is
simple. Let $\rho_s$ be the induced distance on
$\Sigma_s^{(\infty,l)}$ for each $s$, $\|\cdot\|$ be the distance in
the Euclidean metric $R^4$, it is obvious that,
$$
\|x-y\|\leq \rho_s(x,y),~{\rm for}~x,y\in \Sigma_s^{(\infty,l)}.
$$
If there is $C>0$ independent of $s$ such that
$$
\rho_s(x,y)\leq C\|x-y\|,~{\rm for}~x,y\in \Sigma_s^{(\infty,l)},
$$
then $\Sigma_s^{(\infty,l)}$ must be simple, that is,
$\Sigma_{-\infty}^{(\infty,l)}$ is connected.

We conjecture that, {\it every connected component of a blow-up flow
of a symplectic mean curvature flow is simple.}

\section{Blow up analysis at infinity}

In this section we assume that the symplectic mean curvature flow
exists for long time, then we study the structure of the
singularity at infinity.

Suppose that $X_0$ is a blow up point at infinite, then for
arbitrary sequences $t_k\to\infty$ the quantity
\begin{equation}\label{6.1}\max_{\sigma\in (0, r/2]}\sigma^2 \sup_{[t_k-(r-\sigma)^2,
t_k]}\sup_{\Sigma_t\cap B_{r-\sigma}(X_0)}|A|^2\to\infty,
\end{equation}
where $r$ is a constant which is less than injective radius of $M$
at $X_0$. In fact,
\begin{eqnarray*}
\max_{\sigma\in (0, r/2]}\sigma^2 \sup_{[t_k-(r-\sigma)^2,
t_k]}\sup_{\Sigma_t\cap B_{r-\sigma}(X_0)}|A|^2\geq
(r/2)^2\sup_{\Sigma_{t_k}\cap B_{r/2}(X_0)}|A|^2.
\end{eqnarray*}
It is clear that the right hand side term of the above inequality
tends to infinity as $t_k\to \infty$.

We choose $\sigma_k\in (0, r/2]$ such that
$$\sigma_k^2 \sup_{[t_k-(r-\sigma_k)^2, t_k]}
\sup_{\Sigma_t\cap B_{r-\sigma_k}(X_0)}|A|^2=\max_{\sigma\in (0,
r/2]}\sigma^2 \sup_{[t_k-(r-\sigma)^2, t_k]} \sup_{\Sigma_t\cap
B_{r-\sigma}(X_0)}|A|^2.
$$
Let $\tilde t_k\in [t_k-(r-\sigma_k)^2, t_k]$, $F(x_k,\tilde
t_k)=X_k\in\bar{B}_{r-\sigma_k}(X_0)$, such that
$$\lambda_k^2=|A|^2(X_k)=|A|^2(x_k, \tilde t_k)=\sup_{[t_k-(r-\sigma_k)^2, t_k]}
\sup_{\Sigma_t\cap B_{r-\sigma_k}(X_0)}|A|^2.$$ From equation
(\ref{6.1}) we know that $\lambda_k^2\sigma_k^2\to\infty$ as
$k\to\infty$. In particular,
$$\sup_{[t_k-(r-\sigma_k/2)^2, t_k]} \sup_{\Sigma_t\cap
B_{r-\sigma_k/2}(X_0)}|A|^2\leq 4\lambda_k^2,$$ and hence
$$\sup_{[\tilde t_k-(\sigma_k/2)^2, \tilde t_k]}
\sup_{\Sigma_t\cap B_{r-\sigma_k/2}(X_0)}|A|^2\leq 4\lambda_k^2.
$$
Therefore we can consider the rescaled sequence,
$$ F_k(x, s)=\lambda_k(F(x_k+x, \tilde t_k+\lambda_k^{-2}s)-F(x_k, \tilde
t_k)), ~~~~~~~~~~~s\in [-\lambda_k^2\sigma_k^2/4, 0].
$$
We denote the rescaled surface by $\Sigma^k_s=F_k(\Sigma, s)$. By
the same analysis as the one used at a finite time singularity
(see \cite{HL1}), we can show that $\Sigma^k_s\to\Sigma_s^\infty$
in $C^2(B_R(0)\times [-R, R])$ for any $R>0$ and $\Sigma_s^\infty$
is a mean curvature flow in ${\bf C}^2$ which we call a blow-up
flow. It is clear that the blow-up flow also satisfies the
identity and inequalities (\ref{abound}) and (\ref{e0.0}).

Similarly, we can prove the following theorem.

\begin{theorem}\label{infinitetheorem} Let $M$ be a K\"ahler-Einstein
surface and $\Sigma_0$ be a symplectic surface in $M$. Assume that
the mean curvature flow exists globally. Each simple connected
component $\Sigma_{s}^{(\infty,l)}$ of the blow-up flow
$\Sigma^\infty_s$ at infinity is independent of $s$ denoted by
$\Sigma^{(\infty,l)}$, and $\Sigma^{(\infty,l)}$ is a holomorphic
curve in ${\bf C}^2$ with Gauss curvature $-2\leq K\leq 0$, and
finite total curvature
$$
-\int_{\Sigma^\infty}KdV=2\pi N,
$$
where $N$ is a nature number.
\end{theorem}

\section{Singularities of Lagrangian mean curvature flows}

In the case of Lagrangian mean curvature flow, suppose that $M$ is
a K\"ahler-Einstein surface with scalar curvature $\cal R$ and the
mean curvature form of the initial surface $\Sigma_0$ is exact,
Smoczyk \cite{Sm2} showed that there exists a function $\beta$
such that
\begin{eqnarray*}\label{beta}
d\beta &=& H\nonumber\\
\frac{d\beta}{dt} &=& \Delta\beta+\frac{{\cal R}}{4}\beta
\end{eqnarray*}
and $\beta$ is called the Lagrangian angle. Let
$\tilde{\beta}=e^{-\frac{{\cal R}}{4}t}\beta$, we have,
$$\frac{\partial\cos\tilde\beta}{\partial
t}=\Delta\cos\tilde\beta+e^{-\frac{{\cal
R}}{2}t}|H|^2\cos\tilde\beta.
$$
It is certainly more interesting when ${\cal R}=0$, that is, $M$
is a Calabi-Yau surface.

Let $M$ be a compact Calabi-Yau surface with a parallel
holomorphic $(2,0)$-form $\Omega$ of unit length. Then we have
\begin{equation}\label{langle}
\Omega |_{\Sigma_t}= e^{i\beta}d\mu_t =\cos\beta d\mu_t+i\sin\beta
d\mu_t.
\end{equation}

If $\cos\beta>0$, we say that $\Sigma_t$ is almost calibrated.
Smoczyk \cite{Sm2} (also see \cite{TY}) showed that, if $\Sigma_0$
is almost calibrated, then $\Sigma_t$ is too, which we call an
almost calibrated Lagrangian mean curvature flow. It is proved in
\cite{CL2} and \cite{Wa} that, along an almost calibrated
Lagrangian mean curvature flow, there is no Type I singularity.

Let $T$ be an isolated singular time, that is, the mean curvature
flow exists in $t\in [T-\epsilon, T)$, and $(X_0, T)$ be a blow-up
point. We consider the strong convergence of the rescaled surfaces
$\Sigma^k_s$ in $B_R(0)$ around the singular point $X_0$ (c.f.
\cite{CL1}, \cite{HL1}), we have a blow-up flow $\Sigma^\infty_s$
in ${\mathbb{C}}^2$ with the Euclidean metric, which is defined on
$(-\infty, 0]$ for Type I singularity and defined on $(-\infty,
\infty)$ for Type II singularity. The blow-up flow satisfies the
identity and the inequalities (\ref{abound}) and (\ref{e0.0}).

\begin{theorem} Let $M$ be a compact K\"ahler surface and the initial
surface $\Sigma_0$ be an exact Lagrangian surface with $\beta$
bounded. Each simple connected component $\Sigma_{s}^{(\infty,l)}$
of the blow-up flow $\Sigma^\infty_s$ is independent of $s$
denoted by $\Sigma^{(\infty,l)}$, and $\Sigma^{(\infty,l)}$ is a
minimal surface in ${\bf R}^4$ with Gauss curvature $-2\leq K\leq
0$, and finite total curvature
$$
-\int_{\Sigma^\infty}KdV=2\pi N,
$$
where $N$ is a nature number.
\end{theorem}

{\it Proof.} Recall that, by the blow-up technique, we can get a
sequence of rescaled surfaces $\Sigma^k_s$ around $(X_0, T)$,
which converges to $\Sigma^\infty_s$ strongly in $C^2(B_R(0)\times
[-R, R])$ for any $R>0$ and any $B_R(0)\subset {\bf R}^4$, with
$|A_k(0)|=1$ and $|A_k|^2\leq 4$. Moreover, $\Sigma_s^\infty$
evolves along the mean curvature flow. Since
$$
\left(\frac{\partial}{\partial s}-\Delta \right){\beta} _s^k=
\frac{{\cal R}}{4\lambda_k^2}\beta_s^k,
$$ we can see that, along the mean curvature flow $\Sigma_s^\infty$,
$\beta$ satisfies
\begin{equation} \label{4.1} (\frac{\partial}{\partial s}-\Delta)
\beta= 0.
\end{equation}

In the following, we denote the simple connected component
$\Sigma_s^{(\infty,l)}$ of the blow-up flow simply by
$\Sigma_s^{\infty}$.

One can show the monotonicity formula (c.f. \cite{CL2},
\cite{HL1}, \cite{N1},\cite{N2})
\begin{eqnarray*}\label{mon5}
&&\frac{\partial}{\partial s}\left(\int_{\Sigma^\infty_s}\beta^2
\rho(F_\infty, 0, s, s_0)d\mu^\infty_s\right)\nonumber\\&&\leq
-\left( \int_{\Sigma^\infty_s}\beta^2\rho (F_\infty, 0, s, s_0)
\left|H_\infty+\frac{(F_\infty)^{\perp}}{2(s_0-s)}\right|^2d\mu_s^\infty
\right.\nonumber \\
&&\left.+ \int_{\Sigma^\infty_s}|\nabla \beta|^2\rho
(F_\infty,0,s,s_0) d\mu_s^\infty \right) .
\end{eqnarray*}

From this monotonicity formula, similarly we can get that on
$\Sigma^\infty_{-\infty}$,
$$|\nabla \beta|=0, $$
so, $\beta \equiv constant$ on $\Sigma^\infty_{-\infty}$, it can be
normalized to be zero. By the evolution equation (\ref{4.1}), we see
that $\beta\equiv 0$ on $\Sigma^\infty_{s}$. Since
$H_{\Sigma^\infty_s} =d\beta$, we have that
$H_{\Sigma^\infty_s}\equiv 0$, and consequently, $\Sigma^\infty_s$
does not depend on $s$ and is a minimal surface in ${\mathbb{R}}^4$.
Because $\Sigma_s^k$ converges to $\Sigma^\infty$ in $C^2(B_R(0))$
for any $R>0$ as $k\to\infty$, we see that $\Sigma^\infty$ is a
complete minimal surface. Let $K$ be the Gauss curvature of
$\Sigma^\infty$, similarly we can get that
$$\int_{\Sigma^\infty}|A_\infty|^2 d\mu^\infty\leq C,$$
by Proposition $6.1$ in \cite{HO}, we see that that
$$
-\int_{\Sigma^\infty}KdV=2\pi N, $$ where $N$ is a nature number.
This completes the proof of the theorem.

\hfill Q. E. D.

\begin{corollary}
Let $M$ be a compact K\"ahler surface and the initial surface
$\Sigma_0$ be an exact Lagrangian surface with $\beta$ bounded. If
the blow-up flow $\Sigma^\infty_s$ contains only one connected
component, and its limit at $-\infty$, $\Sigma_{-\infty}^\infty$ is
connected, then
 $\Sigma^\infty_s$ is independent of
$s$ denoted by $\Sigma^\infty$, and $\Sigma^\infty$ is a non trivial
minimal surface in ${\bf R}^4$ with Gauss curvature $-2\leq K\leq
0$, and finite total curvature
$$
-\int_{\Sigma^\infty}KdV=2\pi N,
$$
where $N$ is a nature number.
\end{corollary}

By similar arguments, we can show the following theorem. We omit the
proof here.

\begin{theorem} Let $M$ be a compact K\"ahler surface and the initial
surface $\Sigma_0$ be an exact Lagrangian surface with $\beta$
bounded. Assume that the mean curvature flow has a long time
solution. Each simple connected component $\Sigma_{s}^{(\infty,l)}$
of the blow-up flow $\Sigma^\infty_s$ at infinity is independent of
$s$ denoted by $\Sigma^{(\infty,l)}$, and $\Sigma^{(\infty,l)}$ is a
minimal surface in ${\bf R}^4$ with Gauss curvature $-2\leq K\leq
0$, and finite total curvature
$$
-\int_{\Sigma^\infty}KdV=2\pi N,
$$
where $N$ is a nature number.
\end{theorem}

\end{document}